\documentclass[12pt]{article}
\usepackage{amsfonts}
\usepackage{amsmath}
\usepackage{fullpage}
\usepackage{url}

\newcommand{\RR}{\ensuremath{\mathbb{R}}}
\newcommand{\CC}{\ensuremath{\mathbb{C}}}
\newcommand{\HH}{\ensuremath{\mathbb{H}}}
\newcommand{\OO}{\ensuremath{\mathbb{O}}}
\newcommand{\KK}{\ensuremath{\mathbb{K}}}
\newcommand{\mat}[1]{\mathbf{#1}}
\newcommand{\aaa}{\mat{a}}
\newcommand{\bbb}{\mat{b}}
\newcommand{\AAA}{\mat{A}}
\newcommand{\BBB}{\mat{B}}
\newcommand{\II}{\mat{I}}
\newcommand{\MM}{\mat{M}}
\newcommand{\PP}{\mat{P}}
\newcommand{\QQ}{\mat{Q}}
\newcommand{\XX}{\mat{X}}
\newcommand{\GG}{\mat{\Gamma}}
\renewcommand{\SS}{\mat{\hbox{\boldmath$\sigma$}}}
\renewcommand{\bar}[1]{\overline{#1}}
\renewcommand{\star}[1]{#1^*}
\newcommand{\inda}{a} %index 1
\newcommand{\indb}{b} %index 2
\newcommand{\indc}{c} %index 3
\newcommand{\Cl}{\text{C}\ell}
\newcommand{\SO}{\text{SO}}
\newcommand{\SU}{\text{SU}}
\newcommand{\SL}{\text{SL}}
\newcommand{\SP}{\text{SP}}

\newcommand{\tr}{\text{tr}}

\begin{document}

%**end of header

\title{
{\bfseries\boldmath Division Algebra Representations of $\SO(4,2)$}
}

\author{
	Joshua Kincaid and Tevian Dray \\[-2.5pt]
	\normalsize
	\textit{Department of Mathematics, Oregon State University,
		Corvallis, OR  97331} \\[-2.5pt]
	\normalsize
	{\tt kincajos{\rm @}math.oregonstate.edu},
	{\tt tevian{\rm @}math.oregonstate.edu} \\
}

\date{\normalsize 8 August 2013}

\maketitle

\begin{abstract}
Representations of $\SO(4,2)$ are constructed using $4\times4$ and $2\times2$
matrices with elements in $\HH'\otimes\CC$, and the known isomorphism between
the conformal group and $\SO(4,2)$ is written explicitly in terms of the
$4\times4$ representation.
\end{abstract}

% keywords:
% division algebras, magic squares, orthogonal groups, Clifford algebras

% suggested reviewers:
% Tony Sudbery, Robert Wilson, John Baez, John Huerta, Rafal Ablamowicz

\section{Introduction}
\label{intro}

The (local) correspondence between the Lorentz group $\SO(3,1)$ and
$\SL(2,\CC)$ is well-known, and is a natural generalization of the
correspondence between the rotation group $\SO(3)$ and the unitary group
$\SU(2)$.  Manogue and Schray~\cite{Lorentz} generalized this correspondence
to the other division algebras, in particular providing an explicit
construction of $\SO(9,1)$ in terms of $\SL(2,\OO)$.  In subsequent
work~\cite{Dim,Spin}, Manogue and Dray outlined the implications of this
mathematical description for the description of fundamental particles.

Here, we generalize this construction in a different direction, showing how to
describe the conformal group $\SO(4,2)$ as a matrix group over
$\HH'\otimes\CC$, along the way reinterpreting $\SO(3,1)$ as a matrix group
over $\CC'\otimes\CC$.  The resulting parameterization of certain orthogonal
groups in terms of two division algebras is reminiscent of the
Freudenthal-Tits magic square of Lie groups~\cite{Freudenthal,Tits}, and our
results should generalize to the corresponding ``$2\times2$'' version of the
magic square~\cite{SudberyBarton}, shown in Figure~\ref{2x2sq}.

\begin{table}
\small
\begin{center}
\begin{tabular}{|c|c|c|c|c|}
\hline
&$\RR$&$\CC$&$\HH$&$\OO$\\\hline
$\RR'$&
 $\SO(2)\equiv\SU(2,\RR)$&$\SO(3)\equiv\SU(2,\CC)$&
 $\SO(5)\equiv\SU(2,\HH)$&$SO(9)\equiv\SU(2,\OO)$\\
\hline
$\CC'$&
 $\SO(2,1)\equiv\SL(2,\RR)$&$\SO(3,1)\equiv\SL(2,\CC)$&
 $\SO(5,1)\equiv\SL(2,\HH)$&$\SO(9,1)\equiv\SL(2,\OO)$\\
\hline
$\HH'$&
 $\SO(3,2)\equiv\SP(3,\RR)$&$\SO(4,2)\equiv\SU(2,2,\CC)$&
 $\SO(6,2)$&$\SO(10,2)$\\
\hline
$\OO'$&
 $\SO(5,4)$&$\SO(6,4)$&$\SO(8,4)$&$\SO(12,4)$\\
\hline
\end{tabular}
\end{center}
\caption{The ``half-split'' $2\times2$ Lie group magic square.
(The given equivalences are local, that is, up to double-cover.)}
\label{2x2sq}
\end{table}

\section{Split Quaternions}
\label{splitq}

As can be seen from Table~\ref{2x2sq}, the group $\SO(4,2)$ is labeled by the
complex numbers $\CC$ and the split quaternions $\HH'$.  For compatibility
with the notation in~\cite{Dim,Spin}, we write complex numbers $a\in\CC$ in
the form
\begin{equation}
a = x + y \ell
\end{equation}
with $x,y\in\RR$, and we introduce the notation
\begin{equation}
A = z + q K + p KL + t L
\end{equation}
with $p,q,t,z\in\RR$ for split quaternions $A\in\HH'$.  We have of course
\begin{equation}
\ell^2 = -1
\end{equation}
and the split quaternionic multiplication table
\begin{align}
K^2 &= -1; \qquad L^2 = 1 = KL^2 ; \nonumber\\
(K)(L) &= KL; \qquad L(KL) = -K; \qquad (KL)K = L
\end{align}
together with the usual anticommutativity of these unit elements.  The split
quaternions are associative, but not a division algebra, since for
instance
\begin{equation}
(1+L)(1-L) = 0
\end{equation}
so that there are zero divisors.  We introduce separate conjugation operations
\begin{align}
\bar{a} &= x - y \ell \\
\star{A} &= z - q K - p KL - t L
\end{align}

\section{\boldmath The Clifford Algebra $\Cl(4,2)$}
\label{Clifford}

We consider matrices of the form
\begin{equation}
\label{defx}
\XX
  = \begin{pmatrix}A&\bar{a}\\a&-\star{A}\end{pmatrix}
  = \begin{pmatrix}z+qK+pKL+tL&x-y\ell\\
		x+y\ell&-z+qK+pKL+tL\end{pmatrix}
\end{equation}
with $A\in\HH'$ and $a\in\CC$.  Then $\XX$ can be written as
\begin{equation}
\XX = x^\inda\SS_\inda
\label{Xdef}
\end{equation}
where
\begin{equation}
\{x^\inda\} = \{x,y,z,t,p,q\}
\end{equation}
and there is an implicit sum over the repeated index $\inda$.
Equation~(\ref{Xdef}) defines the \textit{generalized Pauli matrices}
\begin{align}
\SS_z &= \begin{pmatrix}1&0\\0&-1\end{pmatrix} ,
  \qquad \SS_x = \begin{pmatrix}0&1\\1&0\end{pmatrix} ,\nonumber\\
\SS_y &= \begin{pmatrix}0&-\ell\\\ell&0\end{pmatrix} ,
  \qquad \SS_t = \begin{pmatrix}L&0\\0&L\end{pmatrix} ,\nonumber\\
\SS_q &= \begin{pmatrix}K&0\\0&K\end{pmatrix} ,
  \qquad\SS_p = \begin{pmatrix}KL&0\\0&KL\end{pmatrix}
\end{align}
which are given this name because $\SS_x$, $\SS_y$, and
$\SS_z$ are just the usual Pauli spin matrices.

We now consider the matrix
\begin{equation}
\PP
  = \begin{pmatrix}0&\XX\\\widetilde{\XX}&0\end{pmatrix}
  = x^\inda\GG_\inda,
\label{Pdef}
\end{equation}
where tilde represents trace reversal,
\begin{equation}
\widetilde{\XX} = \XX - \text{tr}(\XX)\,\II,
\end{equation}
and where the gamma matrices $\GG_\inda$ are implicitly defined
by~(\ref{Pdef}). Explicitly,
\begin{equation}
\GG_\inda
  = \begin{cases}
	\SS_x\otimes\SS_\inda&m\in\{x,y,z\}\\
	\ell\SS_y\otimes\SS_\inda&m\in\{p,q,t\}
    \end{cases}
\end{equation}
or, in more traditional notation,
\begin{align}
\label{4x4exp}
\GG_z
  &= \begin{pmatrix}0&0&1&0\\0&0&0&-1\\1&0&0&0\\0&-1&0&0\end{pmatrix},
 \qquad \GG_x
   = \begin{pmatrix}0&0&0&1\\0&0&1&0\\0&1&0&0\\1&0&0&0\end{pmatrix},
\nonumber\\
\GG_y
  &= \begin{pmatrix}
	0&0&0&-\ell\\0&0&\ell&0\\0&-\ell&0&0\\\ell&0&0&0
     \end{pmatrix},
 \qquad \GG_t
   = \begin{pmatrix}0&0&L&0\\0&0&0&L\\-L&0&0&0\\0&-L&0&0\end{pmatrix},
\nonumber\\
\GG_q
  &= \begin{pmatrix}0&0&K&0\\0&0&0&K\\-K&0&0&0\\0&-K&0&0\end{pmatrix},
 \qquad \GG_p
   = \begin{pmatrix}0&0&KL&0\\0&0&0&KL\\-KL&0&0&0\\0&-KL&0&0\end{pmatrix}.
\end{align}
A straightforward computation using the commutativity of $\CC$ with $\HH'$ now
shows that
\begin{equation}
\{\GG_\inda,\GG_\indb\} = 2g_{\inda\indb}\II
\label{CliffordID}
\end{equation}
where $\II$ is the identity matrix and
\begin{equation}
g_{\inda\indb}
  = \begin{cases}
	0&\inda\neq\indb\\
	1&\inda=\indb\in\{x,y,z,q\}\\
	-1&\inda=\indb\in\{p,t\}
    \end{cases}
\end{equation}
These are precisely the anticommutation relations necessary to generate a
representation of the (real) Clifford algebra $\Cl(4,2)$, so $\PP$
represents an arbitrary element of the vector space underlying $\Cl(4,2)$.

\section{\boldmath Division Algebra Representations of $\SO(4,2)$}
\label{4x4rep}

It is now straightforward to use our representation of the Clifford algebra
$\Cl(4,2)$ to construct a representation of $\SO(4,2)$.  The homogenous
quadratic elements of $\Cl(4,2)$ act as generators of $\SO(4,2)$ via the map
\begin{equation}
\label{paction}
\PP \longmapsto \MM_{\inda\indb}\PP\MM_{\inda\indb}^{-1}
\end{equation}
where
\begin{equation}
\MM_{\inda\indb}
  = \exp\left(-\GG_\inda\GG_\indb\>\frac{\theta}{2}\right)
\end{equation}
and $\PP = x^\inda\GG_\inda$ as above.

In the following, we assume $\inda,\indb,\indc$ are distinct but otherwise
arbitrary indices.  The properties below follow from the Clifford algebra
anticommutation relation~(\ref{CliffordID}):
\begin{align}
\GG_\inda\GG_\inda &= \pm \II,
  \label{prop1}\\
(\GG_\inda\GG_\indb)\GG_\indc
  &= \GG_\indc(\GG_\inda\GG_\indb),
  \label{prop2}\\
(\GG_\inda\GG_\indb)\GG_\indb
  &= (\GG_\indb)^2\GG_\inda
   = g_{\indb\indb}\GG_\inda,
  \label{prop3}\\
(\GG_\inda\GG_\indb)\GG_\inda
  &= -(\GG_\inda)^2\GG_\indb
   = -g_{\inda\inda}\GG_\indb,
  \label{prop4}\\
(\GG_\inda\GG_\indb)^2
  &= -\GG_\inda^2\GG_\indb^2
   = \pm \II
  \label{prop5}.
\end{align}
With these observations, we are prepared to see how $\SO(4,2)$ is generated by
the matrices $\{\GG_\inda\}$.  We compute
\begin{equation}
\label{action4}
\MM_{\inda\indb}\PP\MM_{\inda\indb}^{-1}
  = \exp\left(-\GG_\inda\GG_\indb\>\frac{\theta}{2}\right)
	\left(x^\indc\GG_\indc\right)
	\exp\left(\GG_\inda\GG_\indb\>\frac{\theta}{2}\right).
\end{equation}
From~(\ref{prop1}), if $\inda=\indb$, then $\MM_{\inda\indb}$ is a real
multiple of the identity matrix, which therefore leaves $\PP$ unchanged
under the action~(\ref{paction}).  On the other hand, if $\inda\neq\indb$,
properties~(\ref{prop2})--(\ref{prop4}) imply that $\MM_{\inda\indb}$
commutes with all but two of the matrices $\GG_\indc$.  We therefore
have
\begin{equation}
\GG_\indc\MM_{\inda\indb}^{-1}
  = \begin{cases}
	\MM_{\inda\indb}\GG_\indc,
		&\indc=\inda\text{ or }\indc=\indb\\
	\MM^{-1}_{\inda\indb}\GG_\indc,
		&\inda\neq\indc\neq\indb
    \end{cases}
\label{maction}
\end{equation}
so that the action of $\MM_{\inda\indb}$ on $\PP$ affects only the
$\inda\indb$ plane.  To see what that action is, we first note that if
$\AAA^2 = \pm\II$ then
\begin{equation}
\label{euler}
\exp\left(\AAA\alpha\right)
  = \II\,c(\alpha) + \AAA\,s(\alpha)
  = \begin{cases}
	\II\,\cosh(\alpha) + \AAA\,\sinh(\alpha),&\AAA^2 = \II\\
	\II\,\cos(\alpha) + \AAA\,\sin(\alpha),&\AAA^2 = -\II
    \end{cases}
\end{equation}
where the second equality serves to define the functions $c$ and $s$.
Inserting~(\ref{maction}) and~(\ref{euler}) into~(\ref{action4}), we obtain
\begin{align}
\label{action}
\MM_{\inda\indb} \left(
	x^\inda\GG_\inda + x^\indb\GG_\indb
	\right) \MM_{\inda\indb}^{-1}
  &= \left(\MM_{\inda\indb}\right)^2
	\left(x^\inda\GG_\inda + x^\indb\GG_\indb\right)
  \nonumber\\
  &= \exp\left(-\GG_\inda\GG_\indb\theta\right)
	\left(x^\inda\GG_\inda + x^\indb\GG_\indb\right)
  \nonumber\\
  &= \bigl(
	\II\,c(\theta) - \GG_\inda\GG_\indb\,s(\theta)
     \bigr)
	\left(x^\inda\GG_\inda + x^\indb\GG_\indb\right)
  \nonumber\\
  &= \left(x^\inda c(\theta) - x^\indb s(\theta)g_{\indb\indb}\right)
	\GG_\inda
	+ \left(x^\indb c(\theta) + x^\inda s(\theta)g_{\inda\inda}\right)
	  \GG_\indb.
\end{align}
Thus, the action~(\ref{paction}) is either a rotation from $\inda$ to $\indb$
or a boost in the $\inda\indb$-plane, depending on whether
\begin{equation}
(\GG_\inda\GG_\indb)^2 = \pm\II
\end{equation}
that is, on which version of~(\ref{euler}) is required.   It follows
from~(\ref{Pdef}) that
\begin{align}
\tr(\PP) &= 0 \\
\PP^2 &= - (\det\XX) \, \II
\end{align}
so that the characteristic equation for $\PP$ implies that
\begin{equation}
\det\PP = (\det\XX)^2 = (x^2+y^2+z^2+p^2-q^2-t^2)^2
\end{equation}
which can also be verified by direct computation.  Since transformations of
the form~(\ref{paction}) preserve the determinant of $\PP$, it is clear that
we have constructed $\SO(4,2)$.  The fifteen independent group generators
$\MM_{ab}$ are given explicitly in the appendix.
%in~(\ref{so42mat}).

So far we have considered transformations of the form~(\ref{paction}) acting
on $\PP$, but we can also consider the effect~(\ref{paction}) has on $\XX$.
First, we observe that trace-reversal of $\XX$ corresponds to conjugation in
$\HH'$, that is,
\begin{equation}
\widetilde{\SS_\inda} = \star{\SS}_\inda.
\end{equation}
In light of the off-diagonal structure of the matrices $\GG_\inda$, the
matrices $\GG_\inda\GG_\indb$ then take the block diagonal form
\begin{equation}
\GG_\inda\GG_\indb
  = \begin{pmatrix}
	\SS_\inda\star{\SS}_\indb&0\\
	0&\star{\SS}_\inda\SS_\indb
    \end{pmatrix}
\end{equation}
and, in particular,
\begin{equation}
\exp\left(\GG_\inda\GG_\indb\>\frac{\theta}{2}\right)
  = \begin{pmatrix}
	\exp\left( \SS_\inda\star{\SS}_\indb\>\frac{\theta}{2} \right)&0\\
\noalign{\smallskip}
	0&\exp\left( \star{\SS}_\inda\SS_\indb\>\frac{\theta}{2} \right)
    \end{pmatrix},
\end{equation}
so we can write
\begin{align}
&\exp\left(-\GG_\inda\GG_\indb\>\frac{\theta}{2}\right)
	\>\PP\>
	\exp\left(\GG_\inda\GG_\indb\>\frac{\theta}{2}\right)
	\nonumber\\
&\qquad\qquad=
    \begin{pmatrix}
	0&\exp\left( -\SS_\inda\star{\SS}_\indb\>\frac{\theta}{2} \right) \XX
	  \exp\left( \star{\SS}_\inda\SS_\indb\>\frac{\theta}{2} \right)\\
\noalign{\smallskip}
	\exp\left( -\star{\SS}_\inda\SS_\indb\>\frac{\theta}{2} \right)
		\widetilde{\XX}
	\exp\left(\SS_\inda\star{\SS}_\indb\>\frac{\theta}{2} \right)&0
    \end{pmatrix}.
\end{align}
The $4\times4$ action~(\ref{paction}) acting on $P$ is thus equivalent to the
$2\times2$ action
\begin{equation}
\XX \longmapsto
	\exp\left( -\SS_\inda\star{\SS}_\indb\>\frac{\theta}{2} \right) \XX
	\exp\left( \star{\SS}_\inda\SS_\indb\>\frac{\theta}{2} \right).
\label{xaction}
\end{equation}
on $\XX$.  However, these transformations do not appear to have the
general form
\begin{equation}
\XX\longmapsto\MM\XX\MM^\dagger,
\end{equation}
even if we restrict the dagger operation to include conjugation in just one of
$\HH'$ or $\CC$.
Nonetheless, since $\XX$ is Hermitian with respect to $\CC$, and since
that condition is preserved by~(\ref{xaction}), we will refer to our
$2\times2$ representation of $\SO(4,2)$ as $\SU(2,\HH'\otimes\CC)$.

\section{\boldmath A Real Representation of $SO(4,2)$}
\label{relrep}

We seek now to identify a real representation of $\SO(4,2)$ that satisfies
certain ``nice'' conditions.  We would like our representation to contain a
representation of $\SO(3,1)$ in an obvious way and be linked explicitly to
$\SU(2,\HH'\otimes\CC)$.  Ideally, the construction developed here will extend
naturally to the other groups in Table~\ref{2x2sq} and admit an analog that
can be applied to the Freudenthal-Tits magic square.

We are seeking a real representation, so we require a way to express these as
real matrices while retaining the essential anticommutation relations.  The
solution is provided by finding suitable representations for $\CC$ and $\HH'$
in the form of $2\times2$ real matrices.  We can do this by making use of the
Pauli spin matrices, using the facts that
\begin{align}
(\ell\SS_y)^2 &= -\II,\\
\SS_x^2 &= \SS_z^2 = \II,
\end{align}
and that all three of these matrices are real.  If we map
\begin{align}
1&\to\II\nonumber\\
\ell&\to \ell\SS_y
\end{align}
for $\{1,\ell\}\subset\CC$ and
\begin{align}
1&\to\II\nonumber\\
L &\to \SS_z\nonumber\\
K &\to -\ell\SS_y\nonumber\\
KL &\to\SS_x,
\end{align}
for $\{1, L, K, KL\}\subset\HH'$, then the appropriate multiplication tables
are preserved.  Thus, we can write elements of $\HH'\otimes\CC$ by taking
tensor products of these representations.  We can now write the matrices
$\{\GG_\inda\}$ as
\begin{align*}
\GG_x &= \SS_x\otimes\SS_x\otimes\II\otimes\II,
 \qquad \GG_y = \SS_x\otimes -\ell\SS_y\otimes\II\otimes \ell\SS_y,
 \qquad \GG_z = \SS_x\otimes\SS_z\otimes\II\otimes\II,\\
\GG_t &= \ell\SS_y\otimes\II\otimes\SS_z\otimes\II,
 \qquad \GG_q = \ell\SS_y\otimes\II\otimes-\ell\SS_y\otimes\II,
 \qquad \GG_p = \ell\SS_y\otimes\II\otimes\SS_x\otimes\II.
\end{align*}
In this form, the fifteen generators of $\SO(4,2)$ are
\begin{align}
\label{so42gen}
\GG_t\GG_x &= \SS_z\otimes\SS_x\otimes\SS_z\otimes\II,
 \quad \GG_t\GG_y = -\SS_z\otimes \ell\SS_y\otimes \SS_z\otimes \ell\SS_y,
 \quad \GG_t\GG_z = \SS_z\otimes\SS_z\otimes\SS_z\otimes\II,
 \nonumber\\
\GG_x\GG_y &= \II\otimes\SS_z\otimes\II\otimes \ell\SS_y,
 \quad \GG_y\GG_z = \II\otimes\SS_x\otimes\II\otimes \ell\SS_y,
 \quad \GG_z\GG_x = \II\otimes \ell\SS_y\otimes\II\otimes\II,
 \nonumber\\
\GG_q\GG_x &= -\SS_z\otimes\SS_x\otimes\ell\SS_y\otimes\II,
 \> \GG_q\GG_y = \SS_z\otimes \ell\SS_y\otimes \ell\SS_y \otimes \ell\SS_y,
 \> \GG_q\GG_z = -\SS_z\otimes\SS_z\otimes\ell\SS_y\otimes\II,
 \nonumber\\
\GG_p\GG_x &= \SS_z\otimes\SS_x\otimes\SS_x\otimes\II,
 \quad \GG_p\GG_y = -\SS_z\otimes \ell\SS_y\otimes\SS_x \otimes \ell\SS_y,
 \quad \GG_p\GG_z = \SS_z\otimes\SS_z\otimes\SS_x\otimes\II,
 \nonumber\\
\GG_t\GG_p &= -\II\otimes\II\otimes \ell\SS_y\otimes\II,
 \quad \GG_t\GG_q = \II\otimes\II\otimes\SS_x\otimes\II,
 \quad \GG_p\GG_q = -\II\otimes\II\otimes\SS_z\otimes\II.
\end{align}

There are two steps in this construction at which we can project onto a
representation of $\SO(3,1)$ by projecting from $\HH'$ to $\CC'$.  First, in
our definition of $\XX$ we can set $p = q = 0$, which is equivalent to
restricting $A$ to be in $\CC'\subset\HH'$.  Calling this projection $\pi$, we
then have
\[
\pi(\GG_p) = \pi(\GG_q) = 0
\]
and we are left with only the $t,x,y,$ and $z$ elements.  The anticommutation
relations are obviously still satisfied, and the remaining matrices
$\{\GG_\inda\}$ still generate $\SO(3,1)$ in the obvious way.

On the other hand, we can make the projection from $\HH'$ to $\CC'$ in the
final step by restricting to elements where the third factor is $1 = \II$
or $L = \SS_x$.  However, in this case we get an extra generator,
namely $\GG_{p}\GG_q$.  Exponentiating
$\GG_{p}\GG_q$ to find the corresponding group element
$\MM_{pq}$, we get
\begin{equation}
\MM_{pq} = e^{-L\theta/2}\II,
\end{equation}
which clearly commutes with $\MM_{\inda\indb}\in\SO(3,1)$, so in fact this
is a projection onto
\begin{equation}
\SO(3,1)\times\RR\subset\SO(4,2)
\end{equation}

\section{The Conformal Group}
\label{conformal}

We now show explicitly how to interpret $\SO(4,2)$ as the conformal group by
transforming the representation constructed in Section~\ref{4x4rep} into one
in which the conformal operations are explicit.  The conformal group expands
the Lorentz group, which consists of rotations and boosts, by adding
translations, conformal translations (translations after inverting through the
unit sphere; see~(\ref{conftrans}) below), and a dilation (rescaling).  We
address each of these types of transformations in turn.

Let $V = \text{span}(\{\GG_\inda\})$ and consider $\PP\in V$ as
defined in~(\ref{Pdef}).  We also impose the constraints $p+q\neq0$ and
\begin{equation}
\label{nullcond}
\left|\PP\right|^2
  = \langle\PP,\PP\rangle
  = -t^2 + x^2 + z^2 + y^2 - p^2 + q^2
  = 0.
\end{equation}
where we have introduced the inner product
\begin{equation}
\langle\AAA,\BBB\rangle
% = \frac14 \>\tr(\AAA\circ\BBB)
  = \frac18 \>\tr(\AAA\BBB+\BBB\AAA)
\end{equation}
We then define
\begin{equation}
\QQ = T\GG_t + X\GG_x + Y\GG_y + Z\GG_z,
\end{equation}
with
\begin{align}
\label{Qcoeff}
T
  &= \frac{t}{p + q}
   = -\frac{\langle\GG_t, \PP\rangle}{\langle\GG_p + \GG_q,\PP\rangle},
	\nonumber\\
X
  &= \frac{x}{p + q}
   = \frac{\langle\GG_x, \PP\rangle}{\langle\GG_p + \GG_q,\PP\rangle},
	\nonumber\\
Y
  &= \frac{y}{p + q}
   =\frac{\langle\GG_y, \PP\rangle}{\langle\GG_p + \GG_q,\PP\rangle},
	\nonumber\\
Z
  &= \frac{z}{p + q}
   = \frac{\langle\GG_z, \PP\rangle}{\langle\GG_p + \GG_q,\PP\rangle}.
\end{align}
so that
\begin{equation}
\PP = \QQ(p+q) + p\GG_p + q\GG_q.
\label{PQ}
\end{equation}

We now consider how $\QQ$ changes when elements of $\SO(4,2)$ act on
$\PP$.  As a first observation, when the rotations ($\MM_{xy}$,
$\MM_{yz}$, and $\MM_{zx}$) and boosts ($\MM_{tx}$,
$\MM_{ty}$, and $\MM_{tz}$) act on $\PP$, the effect on $\QQ$
is the same, since $p+q$ is unaffected.
The effect of $\MM_{pq}$ on $p + q$ is given by
\begin{equation}
p + q
  \longmapsto
	p\cosh{\theta} + q\sinh{\theta} + q\cosh{\theta} + p\sinh{\theta}
  = (p + q)(\cosh{\theta} + \sinh{\theta}),
\end{equation}
so that
\begin{equation}
\QQ
  \longmapsto \QQ/(\cosh{\theta} + \sinh{\theta})
  = \QQ e^{-\theta}.
\end{equation}
since $x$, $y$, $z$, and $t$ are unaffected.  This rescaling of $\QQ$ by
$\MM_{pq}$ represents the dilation.

The translations and conformal translations are best understood by considering
null rotations generated by
\begin{equation}
\label{transgen}
\aaa_\inda = \GG_p\GG_\inda - \GG_q\GG_\inda
\end{equation}
and
\begin{equation}
\label{conftransgen}
\bbb_\inda = \GG_p\GG_\inda + \GG_q\GG_\inda
\end{equation}
First, observe that
\begin{equation}
\left(\GG_p\GG_\inda \pm \GG_q\GG_\inda\right)^2
  = \left(\GG_p\GG_\inda\right)^2 + \left(\GG_q\GG_\inda\right)^2
	\pm \GG_p\GG_\inda\GG_p \GG_\inda\GG_q\GG_\inda
	\pm \GG_q\GG_\inda\GG_p\GG_\inda
  = 0,
\end{equation}
where in the last equality we have employed the anticommutation
relations~(\ref{CliffordID}).  As a result,
\begin{equation}
\exp\left(\pm\aaa_\inda\>\frac{\theta}{2}\right)
  = \II \pm \aaa_\inda\,\frac{\theta}{2},
\end{equation}
and
\begin{equation}
\exp\left(\pm\bbb_\inda\>\frac{\theta}{2}\right)
  = \II \pm \bbb_\inda\,\frac{\theta}{2}.
\end{equation}
Next, we compute, as an example, the action of $\aaa_x$ on $\PP$.  To
begin, observe that $\aaa_x$ involves only $\GG_p$,
$\GG_x$, and $\GG_q$, so that $t$, $y$, and $z$ will be
unaffected.  Thus
\begin{align}
\label{transx}
\exp\left(\aaa_x\>\frac{\theta}{2}\right) \GG_x
	\exp\left(-\aaa_x\>\frac{\theta}{2}\right)
  &= \left(\II + \frac{\theta}{2}\GG_p\GG_x
	- \frac{\theta}{2}\GG_q\GG_x\right) \GG_x 
	\left(\II - \frac{\theta}{2}\GG_p\GG_x
	+ \frac{\theta}{2}\GG_q\GG_x\right)\nonumber\\
  &= \left(\GG_x + \frac{\theta}{2}\GG_p
	- \frac{\theta}{2}\GG_q\right)
	\left(\II - \frac{\theta}{2}\GG_p\GG_x
	+ \frac{\theta}{2}\GG_q\GG_x\right)\nonumber\\
  &= \GG_x + \frac{\theta}{2}\GG_p
	- \frac{\theta}{2}\GG_q + \frac{\theta}{2}\GG_p
	+ \frac{\theta^2}{4}\GG_x\nonumber\\
  & \ \ + \frac{\theta^2}{4}\GG_q\GG_p\GG_x
	- \frac{\theta}{2}\GG_q
	+ \frac{\theta^2}{4}\GG_p\GG_q\GG_x
	-\frac{\theta^2}{4}\GG_x\nonumber\\
&=\GG_x + \theta\GG_p - \theta\GG_q.
\end{align}
A similar calculation shows that
\begin{equation}
\label{transp}
\GG_p
  \longmapsto \theta\GG_x - \frac{\theta^2}{2}\GG_q
	+ \left(1 + \frac{\theta^2}{2}\right)\GG_p,
\end{equation}
and
\begin{equation}
\label{transq}
\GG_q
  \longmapsto \GG_q + \theta\GG_x + \frac{\theta^2}{2}\GG_p.
\end{equation}
Combining~(\ref{transx}),~(\ref{transp}), and~(\ref{transq}), we find
\begin{equation}
\label{transsub}
x\GG_x + p\GG_p + q\GG_q
  \longmapsto \left(x + (p + q)\theta\right)\GG_x
	+ \left(p + p\frac{\theta^2}{2} + q\frac{\theta^2}{2} + x\theta\right)
	\GG_p
	+ \left(q - x\theta - p\frac{\theta^2}{2}\right) \GG_q.
\end{equation}
Applying~(\ref{Qcoeff}) to~(\ref{transsub}),
\begin{equation}
X' = \frac{x}{p + q} + \theta = X + \theta.
\end{equation}
In other words, $\aaa_x$ acting on $\PP$ has the effect of translating
$\QQ$ by $\theta$ in the $\GG_x$ direction.  Similar calculations
show that $\aaa_t$, $\aaa_y$, and $\aaa_z$ yield corresponding
translations.

We now consider the effect of $\bbb_x$ acting on $\PP$.  Proceeding as
for~(\ref{transx}), one finds
\begin{equation}
\label{conftransx}
\exp\left(\bbb_x\>\frac{\theta}{2}\right)\GG_x
	\exp\left(-\bbb_x\>\frac{\theta}{2}\right)
  = \GG_x + \frac{\theta}{2}\GG_p + \frac{\theta}{2}\GG_q,
\end{equation}
\begin{equation}
\label{conftransq}
\exp\left(\bbb_x\>\frac{\theta}{2}\right)\GG_q
	\exp\left(-\bbb_x\>\frac{\theta}{2}\right)
  = -\theta\GG_x - \frac{\theta^2}{2}\GG_p
	+ \left(1 - \frac{\theta^2}{2}\right)\GG_q,
\end{equation}
and,
\begin{equation}
\label{conftransp}
\exp\left(\bbb_x\>\frac{\theta}{2}\right)\GG_p
	\exp\left(-\bbb_x\>\frac{\theta}{2}\right)
  = \theta\GG_x + \left(1 + \frac{\theta}{2}\right)\GG_p
	+ \frac{\theta^2}{2}\GG_q.
\end{equation}
Taken together,~(\ref{conftransx}),~(\ref{conftransq}), and~(\ref{conftransp})
yield
\begin{align}
x\GG_x + p\GG_p + q\GG_q
  \longmapsto& \left(x + (p-q)\theta\right)\GG_x
	+ \left(x\theta + p + p\frac{\theta^2}{2} - q\frac{\theta^2}{2}\right)
	\GG_p
\nonumber\\ &\qquad
	+ \left(x\theta + q - q\frac{\theta^2}{2} + p\frac{\theta^2}{2}\right)
	\GG_q,
\end{align}
from which it follows that
\begin{align}
\label{Qconftrans}
X'
  &= \frac{x + (p-q)\theta}
	{x\theta + p + p\frac{\theta^2}{2} - q\frac{\theta^2}{2}
	+ x\theta + q - q\frac{\theta^2}{2} + p\frac{\theta^2}{2}}\nonumber\\
  &= \frac{X + \frac{p-q}{p+q}\theta}
	{2X\theta + 1 + \frac{p-q}{p+q}\theta^2}\nonumber\\
  &= \frac{X + |\QQ|^2\theta}{2X\theta + 1 + |\QQ|^2\theta^2},
\end{align}
where in the last line we have used~(\ref{nullcond}) to write
\begin{align}
\frac{p - q}{p + q}
  &= \frac{p^2 - q^2}{(p + q)^2}
   = \frac{-t^2 + x^2 + y^2 + z^2}{(p+q)^2}\nonumber\\
  &= -T^2 + X^2 + Y^2 + Z^2\nonumber\\
  &= \langle \QQ,\QQ\rangle \equiv |\QQ|^2.
\end{align}
To see why this is a conformal translation, we note that an element $v$ in an
inner product space $V$ satisfies
\begin{equation}
v^{-1} = \frac{v}{|v|^2}.
\end{equation}
Then a conformal translation of $v$ in the direction of the vector $\alpha$ is
given by
\begin{equation}
\label{conftrans}
v
  \longmapsto \left(v^{-1} + \alpha\right)^{-1}
  = \frac{v + \alpha|v|^2}{1 + 2\langle v,\alpha\rangle + |\alpha|^2|v^2}.
\end{equation}
Taking $v = \QQ$ and assuming $\alpha = \theta\GG_x$, the
$\GG_x$ component of~(\ref{conftrans}) becomes precisely
(\ref{Qconftrans}).  Again, a similar calculation shows that $\bbb_y$,
$\bbb_z$, and $\bbb_t$ are the other conformal translations.

The fifteen elements of $\SO(4,2)$ therefore act on $\QQ$
via~(\ref{paction}) and~(\ref{PQ}) as the conformal group, with 3 rotations, 3
boosts, a dilation, 4 translations, and 4 conformal translations.

\section{Conclusion}
\label{conclude}

As shown in Section~\ref{4x4rep}, we obtain a $2\times2$ representation of
$\SO(4,2)$ over $\HH'\otimes\CC$ simply by restricting the $4\times4$
representation via~(\ref{paction}) to one of the $2\times2$ blocks of
$\PP$, say $\XX$.  The resulting action via~(\ref{xaction}) reproduces
that of Manogue and Schray~\cite{Lorentz} when both are restricted to
$\SO(3,1)$, that is, when the matrix elements are complex.  However, as
already noted, it is not possible to express~(\ref{xaction}) in the simple
form involving Hermitian conjugation that was used by Manogue and Schray in
all cases.  Nonetheless, we have shown that
\begin{equation}
\SO(4,2) \equiv \SU(2,\HH'\otimes\CC)
\end{equation}
(up to double-cover issues).  This is our primary result.

We expect the construction given in Sections~\ref{Clifford} and~\ref{4x4rep}
to carry over with minimal modification to any combination of division
algebras, yielding a representation of each group in Table~\ref{2x2sq} of the
form $SU(2,\KK'\otimes\KK)$.  For further details, see~\cite{2x2}.  It is
hoped that this construction can be further generalized to the
Freudenthal-Tits magic square itself, leading to an interpretation of each
group in that magic square of the form $SU(3,\KK'\otimes\KK)$, and providing
new insight into the exceptional groups $E_6$, $E_7$, and $E_8$.

\section*{Acknowledgments}

This paper is based on a thesis submitted by JK in partial fulfillment of the
degree requirements for his M.S.\ in Mathematics at Oregon State
University~\cite{JoshuaThesis}.  We thank John Huerta and Corinne Manogue for
their comments on early versions of this work.  The completion of this paper
was made possible in part through the support of a grant from the John
Templeton Foundation.

\section*{Appendix}

We list here the fifteen $4\times4$ matrices over $\HH'\otimes\CC$ that
generate $SO(4,2)$ under the action~(\ref{paction}):

{\small
\begin{align}
\MM_{xy}
  &= \left(\begin{smallmatrix}
	e^{i\phi/2}&0&0&0\\
	0&e^{-i\phi/2}&0&0\\
	0&0&e^{i\phi/2}&0\\
	0&0&0&e^{-i\phi/2}
      \end{smallmatrix}\right),\nonumber\displaybreak[0]\\
\MM_{yz}
  &= \left(\begin{smallmatrix}
	\cos{\frac{\phi}{2}}&i\sin{\frac{\phi}{2}}&0&0\\
	i\sin{\frac{\phi}{2}}&\cos{\frac{\phi}{2}}&0&0\\
	0&0&\cos{\frac{\phi}{2}}&i\sin{\frac{\phi}{2}}\\
	0&0&i\sin{\frac{\phi}{2}}&\cos{\frac{\phi}{2}}
      \end{smallmatrix}\right),\nonumber\displaybreak[0]\\
\MM_{zx}
  &= \left(\begin{smallmatrix}
	\cos{\frac{\phi}{2}}&\sin{\frac{\phi}{2}}&0&0\\
	-\sin{\frac{\phi}{2}}&\cos{\frac{\phi}{2}}&0&0\\
	0&0&\cos{\frac{\phi}{2}}&\sin{\frac{\phi}{2}}\\
	0&0&-\sin{\frac{\phi}{2}}&\cos{\frac{\phi}{2}}
      \end{smallmatrix}\right),\nonumber\displaybreak[0]\\
\MM_{qx}
  &= \left(\begin{smallmatrix}
	\cos{\frac{\phi}{2}}&K\sin{\frac{\phi}{2}}&0&0\\
	K\sin{\frac{\phi}{2}}&\cos{\frac{\phi}{2}}&0&0\\
	0&0&\cos{\frac{\phi}{2}}&-K\sin{\frac{\phi}{2}}\\
	0&0&-K\sin{\frac{\phi}{2}}&\cos{\frac{\phi}{2}}
      \end{smallmatrix}\right),\nonumber\displaybreak[0]\\ 
\MM_{qy}
  &= \left(\begin{smallmatrix}
	\cos{\frac{\phi}{2}}&-K\otimes i\sin{\frac{\phi}{2}}&0&0\\
	K\otimes i\sin{\frac{\phi}{2}}&\cos{\frac{\phi}{2}}&0&0\\
	0&0&\cos{\frac{\phi}{2}}&K\otimes i\sin{\frac{\phi}{2}}\\
	0&0&-K\otimes i\sin{\frac{\phi}{2}}&\cos{\frac{\phi}{2}}
      \end{smallmatrix}\right),\nonumber\displaybreak[0]\\
\MM_{qz}
  &= \left(\begin{smallmatrix}
	e^{K\phi/2}&0&0&0\\
	0&e^{-K\phi/2}&0&0\\
	0&0&e^{-K\phi/2}&0\\
	0&0&0&e^{K\phi/2}
      \end{smallmatrix}\right),\nonumber\displaybreak[0]\\
\MM_{tp}
  &= \left(\begin{smallmatrix}
	e^{K\phi/2}&0&0&0\\
	0&e^{K\phi/2}&0&0\\
	0&0&e^{K\phi/2}&0\\
	0&0&0&e^{K\phi/2}
      \end{smallmatrix}\right),\nonumber\displaybreak[0]\\
\MM_{tx}
  &= \left(\begin{smallmatrix}
	\cosh{\frac{\phi}{2}}&L\sinh{\frac{\phi}{2}}&0&0\\
	L\sinh{\frac{\phi}{2}}&\cosh{\frac{\phi}{2}}&0&0\\
	0&0&\cosh{\frac{\phi}{2}}&-L\sinh{\frac{\phi}{2}}\\
	0&0&-L\sinh{\frac{\phi}{2}}&\cosh{\frac{\phi}{2}}
      \end{smallmatrix}\right),\nonumber\displaybreak[0]\\
\MM_{ty}
  &= \left(\begin{smallmatrix}
	\cosh{\frac{\phi}{2}}&-L\otimes i\sinh{\frac{\phi}{2}}&0&0\\
	L\otimes i\sinh{\frac{\phi}{2}}&\cosh{\frac{\phi}{2}}&0&0\\
	0&0&\cosh{\frac{\phi}{2}}&-L\otimes i\sinh{\frac{\phi}{2}}\\
	0&0&L\otimes i\sinh{\frac{\phi}{2}}&\cosh{\frac{\phi}{2}}
      \end{smallmatrix}\right),\nonumber\displaybreak[0]\\
\MM_{tz}
  &= \left(\begin{smallmatrix}
	e^{L\phi/2}&0&0&0\\
	0&e^{-L\phi/2}&0&0\\
	0&0&e^{-L\phi/2}&0\\
	0&0&0&e^{L\phi/2}
      \end{smallmatrix}\right),\nonumber\displaybreak[0]\\
\MM_{tq}
  &= \left(\begin{smallmatrix}
	e^{KL\phi/2}&0&0&0\\
	0&e^{KL\phi/2}&0&0\\
	0&0&e^{KL\phi/2}&0\\
	0&0&0&e^{KL\phi/2}
      \end{smallmatrix}\right),\nonumber\displaybreak[0]\\
\MM_{px}
  &= \left(\begin{smallmatrix}
	\cosh{\frac{\phi}{2}}&KL\sinh{\frac{\phi}{2}}&0&0\\
	KL\sinh{\frac{\phi}{2}}&\cosh{\frac{\phi}{2}}&0&0\\
	0&0&\cosh{\frac{\phi}{2}}&-KL\sinh{\frac{\phi}{2}}\\
	0&0&-KL\sinh{\frac{\phi}{2}}&\cosh{\frac{\phi}{2}}
      \end{smallmatrix}\right),\nonumber\displaybreak[0]\\
\MM_{py}
  &= \left(\begin{smallmatrix}
	\cosh{\frac{\phi}{2}}&-KL\otimes i\sinh{\frac{\phi}{2}}&0&0\\
	KL\otimes i\sinh{\frac{\phi}{2}}&\cosh{\frac{\phi}{2}}&0&0\\
	0&0&\cosh{\frac{\phi}{2}}&KL\otimes i\sinh{\frac{\phi}{2}}\\
	0&0&-KL\otimes i\sinh{\frac{\phi}{2}}&\cosh{\frac{\phi}{2}}
      \end{smallmatrix}\right),\nonumber\displaybreak[0]\\
\MM_{pz}
  &= \left(\begin{smallmatrix}
	e^{KL\phi/2}&0&0&0\\
	0&e^{-KL\phi/2}&0&0\\
	0&0&e^{-KL\phi/2}&0\\
	0&0&0&e^{KL\phi/2}
      \end{smallmatrix}\right),\nonumber\displaybreak[0]\\ 
\MM_{pq}
  &= \left(\begin{smallmatrix}
	e^{-L\phi/2}&0&0&0\\
	0&e^{-L\phi/2}&0&0\\
	0&0&e^{-L\phi/2}&0\\
	0&0&0&e^{-L\phi/2}
      \end{smallmatrix}\right).
\label{so42mat}
\end{align}
}

\newpage
%\bibliographystyle{unsrt}
%\let\TDbib=1\bibliographystyle{TDbib}
%\bibliography{Octonions,octonions}

\end{document}